\documentclass[12pt]{amsart}

\usepackage{amsmath}
\usepackage{amsfonts}

\newcommand{\riem}[1][]{\ensuremath{\Sigma_{#1}}}  
\newcommand{\st}{\, | \,} 
\newcommand{\Oqc}{\mathcal{O}_{\mathrm{qc}}} 

\theoremstyle{plain}
        \newtheorem{theorem}{Theorem}[section]
        \newtheorem{lemma}[theorem]{Lemma}
        \newtheorem{proposition}[theorem]{Proposition}
        \newtheorem{corollary}[theorem]{Corollary}

\theoremstyle{definition}
        \newtheorem{definition}[theorem]{Definition}

\theoremstyle{remark}
    \newtheorem{remark}[theorem]{Remark}

\numberwithin{equation}{section} 

\usepackage{fullpage}

\begin{document}

\title{A complex structure on the set of quasiconformally extendible
non-overlapping mappings into a Riemann surface}

\date{\today}

\author{David Radnell}
\address{Department of Mathematics and Statistics \\
American University of Sharjah \\
PO Box 26666, Sharjah, UAE} \email[D. ~Radnell]{dradnell@aus.edu}

\author{Eric Schippers}
\address{Department of Mathematics \\
University of Manitoba\\
Winnipeg, MB, R3T 2N2, Canada} \email[E.
~Schippers]{eric\_schippers@umanitoba.ca}

\begin{abstract}
 Let $\riem$ be a Riemann surface with $n$ distinguished
 points $p_1,\ldots,p_n$.  We prove that the set of $n$-tuples $(\phi_1,\ldots,\phi_n)$ of
 univalent mappings $\phi_i$ from the unit disc $\mathbb{D}$
 into $\riem$ mapping $0$ to $p_i$, with non-overlapping images and
 quasiconformal extensions to a neighbourhood of
 $\overline{\mathbb{D}}$, carries a natural complex Banach manifold structure.
 This complex structure is locally modelled on the $n$-fold
 product of a two complex-dimensional extension of the universal
 Teichm\"uller space.  Our results are motivated by Teichm\"uller
 theory and two-dimensional conformal field theory.
\end{abstract}


\maketitle



\begin{section}{Introduction}
\label{introduction}

Let $\overline{\mathbb{C}}$ be the Riemann sphere, $\mathbb{D} = \{ z \in \mathbb{C} \,|\, |z|<1\}$ and $\mathbb{D}^*=\{z \in \overline{\mathbb{C}} \, | \, |z|>1\}$.

\begin{definition}[of $\mathcal{O}_{qc}(\riem)$] \label{de:Oqcdefinition}
 Let $\riem$ be a Riemann surface with $n$ distinguished points $p_1,\dotsc,p_n$.
 Define $\Oqc(\riem)$ to be the set of $n$-tuples
 $(\phi_1,\ldots,\phi_n)$
 where $\phi_i \colon\mathbb{D} \rightarrow \riem$ are maps
 with the following
 properties:
 \begin{enumerate}
  \item $\phi_i(0)=p_i$
  \item $\phi_i$ is one-to-one and holomorphic on $\mathbb{D}$
  \item $\phi_i$ has a quasiconformal extension to an open neighbourhood of
   $\overline{\mathbb{D}}$
  \item $\phi_i(\overline{\mathbb{D}}) \cap
  \phi_j(\overline{\mathbb{D}}) = \emptyset$ whenever $i \neq j$.
 \end{enumerate}
\end{definition}

\begin{remark} The space $\mathcal{O}_{qc}(\riem)$ was
 introduced by the authors in \cite[Definition 5.1]{RS05}. In this
 paper the direction of the maps has been reversed and the notation
 slightly modified.
\end{remark}

\begin{remark}  The motivation for this work came from the case of compact $\riem$. The authors thank the referee for carefully reading the manuscript and for pointing out that our results hold when $\riem$ is an arbitrary Riemann surface.
\end{remark}

\begin{definition} The model space $\Oqc$ is defined as follows:
\begin{equation*}
 \begin{split}
\Oqc = \{ f \colon \mathbb{D} \rightarrow \mathbb{C}
   \; | \; & f
 \ \text{is one-to-one, holomorphic, has quasiconformal extension to } \mathbb{C}, \\
 & \text{ and } f(0)=0 \}.
 \end{split}
 \end{equation*}
\end{definition}

Our main result is that $\Oqc(\riem)$ is a complex Banach manifold
with complex structure locally modelled on $\Oqc \times \cdots
\times \Oqc$ (Theorem \ref{th:Oqcriem_complexstructure}).

Non-overlapping univalent mappings are a natural object of study
in geometric function theory and have been extensively researched,
especially in connection with extremal problems and inequalities.
If $\riem = \overline{\mathbb{C}}$ with $n$ distinguished points,
then the class $\Oqc(\riem)$ with the restriction that the maps
have quasiconformal extension removed is referred to as the
Goluzin-Lebedev class.

Non-overlapping maps with quasiconformal
extension have also been studied a great deal, especially in the
case of pairs of non-overlapping maps.  The condition that the maps
have quasiconformal extension usually results in sharpened
inequalities. A recent review of the literature on non-overlapping
maps can be found in \cite{Grinshpan}.

In the case of several variables, Forstneri\v c \cite{Forstneric}
showed that the set of holomorphic mappings of pseudoconvex
domains into a complex manifold possesses a complex structure, for
various choices of regularity of the mappings on the boundary.

In this paper we consider the problem of the existence of a complex
structure on $\Oqc(\riem)$.  For compact $\riem$, this question is motivated by both
Teichm\"uller theory and conformal field theory.  In conformal field
theory, a fundamental object is the rigged moduli space. This infinite-dimensional moduli space is the space of conformal equivalence classes of Riemann surfaces together with specified local coordinates at the $n$ distinguished points (see \cite{RS05} for details and references).  In \cite{RS05} we showed,
by generalizing to elements of $\Oqc(\riem)$ as the local coordinate
data at the distinguished points, that this moduli space has a natural complex
structure inherited from the infinite-dimensional Teichm\"uller
space of Riemann surfaces with $n$ boundary curves. Furthermore, our
results in \cite{RS05} suggested that the Teichm\"uller space of
genus $g$ Riemann surfaces with $n$ boundary curves is fibred over
the (finite-dimensional) Teichm\"uller space of $n$-punctured, genus $g$ surfaces.
The fibres are $\Oqc(\riem)$, so these fibres should
inherit a complex structure.

In this paper we show that $\Oqc(\riem)$ possesses a natural
complex structure arising from a different point of view.  We
construct this complex structure from a two-complex-dimensional
extension of the universal Teichm\"uller space using recent
results of Teo \cite{Teo}. Interestingly, our construction does
not have any apparent connection with the fact that the
Teichm\"uller space of a Riemann surface with $n$ boundary curves
is contained in the universal Teichm\"uller space. Nevertheless,
we suspect that these two complex structures are compatible, and
we hope to settle the question in a future publication.
\end{section}


\begin{section}{Preliminaries}

\label{se:BerspreBers}
 In this section we
summarize some known results about the Bers embedding which will be
needed in the rest of the paper. Some of these results are standard,
but we also rely on results of Teo \cite{Teo}. Details can be found
in \cite{Bers}, \cite{Lehto}, \cite{Nag} and \cite{Teo}.

 Using the normalizations in \cite{Teo} and \cite{TTmem} we define the
 following function
 spaces:
 \begin{equation*}
 \begin{split}
   \mathcal{D}  = \{ f  \colon \mathbb{D} \rightarrow \overline{\mathbb{C}} \; | \; & f \text{ is one-to-one, holomorphic, has quasiconformal extension to } \overline{\mathbb{C}}, \\
   & \text{ and } f(0)=0, f'(0)=1, f''(0)=0 \}
 \end{split}
 \end{equation*}
 and
\begin{equation*}
 \begin{split}
\widetilde{\mathcal{D}} = \{ f \colon \mathbb{D} \rightarrow \mathbb{C}
   \; | \; & f
 \ \text{is one-to-one, holomorphic, has quasiconformal extension to } \mathbb{C}, \\
 & \text{ and } f(0)=0, f'(0)=1 \}.
 \end{split}
 \end{equation*}

 The function space $\mathcal{D}$ is one of several
 standard models of the universal Teichm\"uller space.
 Let $L^\infty_1(\mathbb{D}^*)$ denote the set of measurable
 functions on $\mathbb{D}^*$ with sup norm strictly less than one.
 We define an equivalence relation on $L^\infty_1(\mathbb{D}^*)$ as
 follows: $\mu \sim \nu$ if and only if the normalized solutions
 $w_\mu \colon\mathbb{D}^* \rightarrow \mathbb{D}^*$ and
 $w_\nu \colon\mathbb{D}^* \rightarrow \mathbb{D}^*$ to the Beltrami
 equation satisfy $w_\mu=w_\nu$ on $\partial \mathbb{D}^*$.   We
 then define the universal Teichm\"uller space
 \[  T(\mathbb{D}^*)= L^\infty_1(\mathbb{D}^*)/\sim.  \]
 Furthermore, $\mu \sim \nu$
 if and only if the solutions
 $f_\mu$ and $f_\nu$ to the Beltrami equation on
 $\overline{\mathbb{C}}$ with $\mu$ and $\nu$ extended to be $0$ on $\mathbb{D}$
 (normalized so that $f(0)=f'(0)-1=f''(0)=0$) satisfy
 $f_\mu=f_\nu$ on $\mathbb{D}$.  Thus, since their restrictions to
 $\mathbb{D}$ are holomorphic and one-to-one, $T(\mathbb{D}^*)$ is
 in one-to-one correspondence with $\mathcal{D}$.

\begin{remark}
 The normalization for $\mathcal{D}$ is
 not the usual one, but it is convenient in our setting.
\end{remark}
 \begin{remark}
  In \cite{Teo} and \cite{TTmem}, the space $\mathcal{D}$
  is defined slightly differently: functions in $\mathcal{D}$
  have range confined to $\mathbb{C}$.  This seems to be a minor
  oversight, since if $\mathcal{D}$ is
  to be identified with $T(\mathbb{D}^*)$, then the normalization $f''(0)=0$
  cannot be
  imposed without allowing maps in $\mathcal{D}$ to take the value
  $\infty$.
 \end{remark}

 The space $\widetilde{\mathcal{D}}$ is a model of
 the universal Teichm\"uller curve $\mathcal{T}(\mathbb{D}^*)$ of Bers
  \cite{Bers, Teo, TTmem}.
  Although the identification of $T(\mathbb{D}^*)$ with $\mathcal{D}$ will
  be directly used in this paper, the identification of
  $\mathcal{T}(\mathbb{D}^*)$ with $\widetilde{\mathcal{D}}$
  will not be.  However, some of the arguments in this paper regarding the
  complex structures of the spaces could equally well be phrased in
  terms of $\mathcal{T}(\mathbb{D}^*)$
  rather than $T(\mathbb{D}^*)$.

 \begin{remark}
  By conformal welding one can make the identifications
 \begin{align*}
  \mathcal{D} & \cong  T(\mathbb{D}^*) \cong \mbox{QS}(S^1)/
  \mbox{M\"ob}(S^1) \\
  \intertext{and}
  \widetilde{\mathcal{D}} & \cong  \mathcal{T}(\mathbb{D}^*) \cong
  \mbox{QS}(S^1)/ \mbox{Rot}(S^1),
 \end{align*}
 where $QS(S^1)$ denotes the set of quasisymmetric maps from the
 unit circle $|z|=1$ to itself, and
 $\mbox{Rot}(S^1)$ is the set of transformations $z \mapsto
 e^{i\theta} z$.
 In the Appendix (Section \ref{se:Oqcappendix})  it is
 shown that there is a natural bijection
  \[  \Oqc \cong \mbox{QS}(S^1) \times \mathbb{R}, \]
  which can be seen as the next step in this progression
  of removing normalizations from $T(\mathbb{D}^*)$.
   This fact
  is not used in the rest of the paper, but is an interesting point
  of interpretation.
 \end{remark}

We briefly outline how $\mathcal{D}$ and $\widetilde{\mathcal{D}}$
can be mapped respectively into open subsets of the Banach spaces
\[  A^2_\infty(\mathbb{D}) = \{ \psi(z) \colon \mathbb{D} \rightarrow
 \mathbb{C} \;|\; \psi \mbox{ holomorphic}, \ ||\psi||_{2,\infty}=
 \sup_{z \in \mathbb{D}}
 (1- |z|^2)^2 | \psi(z)| < \infty \}  \]
and
\[  A^1_\infty(\mathbb{D}) = \{ \phi(z) \colon \mathbb{D} \rightarrow
 \mathbb{C} \;|\; \phi \mbox{ holomorphic}, \ ||\phi||_{1,\infty}=
 \sup_{z \in \mathbb{D}}
 (1- |z|^2) | \phi(z)| < \infty \}.  \]

\begin{remark}
 $A^2_\infty(\mathbb{D})$ should be interpreted as a space of
 holomorphic quadratic differentials $\psi(z)dz^2$, and
 $A^1_\infty(\mathbb{D})$
 should be interpreted as a space of first-order holomorphic
 differentials $\phi(z)dz$.
\end{remark}
We define two differential operations on holomorphic functions
\begin{align*}
 \mathcal{A}(f) &= \frac{f''}{f'} \\
 \intertext{and}
 \mathcal{S}(f) &=  \frac{f'''}{f'} - \frac{3}{2}
 \left(\frac{f''}{f'}\right)^2
\end{align*}
which are called the ``pre-Schwarzian derivative'' and the
``Schwarzian derivative'' respectively. The Bers embedding
$\mathcal{S} \colon\mathcal{D} \rightarrow A^2_\infty(\mathbb{D})$ maps
the universal Teichm\"uller space injectively into the Banach space
$A^2_\infty(\mathbb{D})$. A well-known result of Ahlfors says that
the image is open, and thus $T(\mathbb{D}^*)$ and $\mathcal{D}$
inherit a complex structure. This complex structure is compatible
with that inherited from $L^\infty_1(\mathbb{D}^*)$ by results of
Bers \cite{Bers}.

We now define the natural injections of $\widetilde{\mathcal{D}}$
into the Banach spaces $A^1_\infty(\mathbb{D})$ and
$A^2_\infty(\mathbb{D})\oplus \mathbb{C}$.  Firstly, there is the
``Bers embedding'' of the universal Teichm\"uller curve:
\begin{eqnarray*}
 \beta \colon \widetilde{\mathcal{D}} &\rightarrow& A^2_\infty(\mathbb{D})
 \oplus \mathbb{C} \\
 f &\mapsto& \left( \mathcal{S}(f),\frac{1}{2}\mathcal{A}(f)(0)
 \right).
\end{eqnarray*}
Bers showed that the image of $\beta$ is open \cite{Bers}. Secondly,
there is the injection defined by
\begin{eqnarray*}
 \hat{\beta}  \colon \widetilde{\mathcal{D}} & \rightarrow &
 A^1_\infty(\mathbb{D})
 \\
  f &\mapsto& \mathcal{A}(f).
\end{eqnarray*}
 Using these injections we define metrics on
 $\mathcal{D}$ and $\widetilde{\mathcal{D}}$
 as follows.
 For $f_1,f_2 \in \mathcal{D}$, let
 \[  d_s(f_1,f_2)=
 ||\mathcal{S}(f_1)-\mathcal{S}(f_2)||_{2,\infty}.  \]
 For $f_1, f_2 \in \widetilde{\mathcal{D}}$, let
 \[
 d_{ps}(f_1,f_2)=||\mathcal{A}(f_1)-\mathcal{A}(f_2)||_{1,\infty}.
  \]

\begin{remark} \label{re:SofDtildeopen}
 By the invariance of the Schwarzian under
 post-composition by M\"obius transformations,
 $\mathcal{S}(\widetilde{\mathcal{D}})=\mathcal{S}(\mathcal{D})$.  In
 particular $\mathcal{S}(\widetilde{\mathcal{D}})$ is open in
 $A^2_\infty(\mathbb{D})$.
\end{remark}

We require results of Teo \cite{Teo} concerning the compatibility
of the complex structures on $\widetilde{\mathcal{D}}$ induced
from $A^1_\infty(\mathbb{D})$ and $A^2_\infty(\mathbb{D})\oplus
\mathbb{C}$ by $\hat{\beta}$ and $\beta$ respectively. These are
summarized in the following four theorems.

Define maps
\begin{align*}
 \Psi \colon A^1_\infty(\mathbb{D}) &\rightarrow A^2_\infty(\mathbb{D}) \\
 g & \mapsto  g'-\frac{g^2}{2} \\
 \intertext{and}
 \hat{\Psi}  \colon A^1_\infty(\mathbb{D}) & \rightarrow
 A^2_\infty(\mathbb{D}) \oplus \mathbb{C}  \\
 g & \mapsto  \left(\Psi(g),\frac{1}{2}g(0)\right).
 \end{align*}
The importance of these maps stems from the identities $\Psi \circ
\mathcal{A} = \mathcal{S}$ and $\hat{\Psi} \circ \hat{\beta} = \beta$.

\begin{theorem}\cite[Corollary A.3]{Teo} \label{th:Psicontinuous}
 $\Psi$ is holomorphic.
\end{theorem}

\begin{theorem}\cite[Theorem A.5]{Teo}
\label{th:hatpsiholomorphic}
 $\hat{\Psi}$ is one-to-one and holomorphic.
\end{theorem}

\begin{theorem} \label{th:hatbetaopen}
 The map $\hat{\beta}$ is one-to-one and the image
 is open in $A^1_\infty(\mathbb{D})$.  In particular
 $\widetilde{\mathcal{D}}$ inherits a complex structure from
 $A^1_\infty(\mathbb{D})$.
\end{theorem}
\begin{proof}
 $\hat{\beta}$ is clearly one-to-one, since if
 $\mathcal{A}(f)=\mathcal{A}(g)$ we have that in some open set
 \[  \log{f'} = \log{g'} + C \]
 for some branch of logarithm and constant $C$ (note $f'\neq 0$ and $g' \neq 0$
 since they are univalent), so $f'=e^C g'$.
 Since $f'(0)=g'(0)$ and $f(0)=g(0)=0$ it follows that $f=g$.

 The image of $\beta$
 is open in $A^2_\infty(\mathbb{D}) \oplus \mathbb{C}$ by a result
 of Bers \cite{Bers}, as mentioned above.  Since
 $\hat{\Psi} \circ \hat{\beta}=\beta$ and
 $\hat{\Psi}$ is continuous and one-to-one, this implies that the
 image of $\hat{\beta}$ is open.
\end{proof}

\begin{theorem} \label{th:BersequalspreBers}
 The restriction $\hat{\Psi} \colon\hat{\beta}(\widetilde{\mathcal{D}}) \rightarrow
 \beta(\widetilde{\mathcal{D}})$ is a biholomorphism.
\end{theorem}
\begin{proof}  This is stated and proved in the proof of
 \cite[Theorem A.6]{Teo}.
\end{proof}
\begin{remark}  In particular, $A^1_\infty(\mathbb{D})$ and
 $A^2_\infty(\mathbb{D})$ induce the same complex structure on
 $\widetilde{\mathcal{D}}$ via $\hat{\beta}$ and $\beta$
 respectively.
 This is the actual statement of
 \cite[Theorem A.6]{Teo}.
\end{remark}
\begin{remark} Teo's Theorem \ref{th:BersequalspreBers} is not an
 immediate
 consequence of Theorems \ref{th:hatpsiholomorphic} and
 \ref{th:hatbetaopen}.  One-to-one holomorphic maps in
 infinite dimensions do not necessarily have holomorphic inverses.
\end{remark}

\end{section}


\begin{section}{Complex structure on $\Oqc$ and $\Oqc(\riem)$}
\begin{subsection}{Complex structure on $\Oqc$}
In the rest of the paper, we will always consider $\mathcal{D}$
and $\widetilde{\mathcal{D}}$ to have the topology induced by the
metrics $d_s$ and $d_{ps}$.

We now define an embedding of $\mathcal{O}_{qc}$ into a Banach
space:
\begin{eqnarray*}
 \chi \colon \mathcal{O}_{qc} & \rightarrow & A^1_\infty(\mathbb{D})
 \oplus \mathbb{C}  \\
 f & \mapsto & \left( \mathcal{A}(f), f'(0) \right).
\end{eqnarray*}
The Banach space direct sum norm on
$A^1_\infty(\mathbb{D})\oplus \mathbb{C}$ is defined by
$
||(\phi, c)||=||\phi||_{1,\infty} +
|c|.
$
\begin{theorem} \label{th:chi_induces_complex_structure}
 $\chi$ is one-to-one, and the image is $\hat{\beta}(\widetilde{\mathcal{D}})
 \oplus \mathbb{C}^*$.  In particular, $A^1_\infty \oplus \mathbb{C}$
 induces a complex Banach manifold structure on $\Oqc$ via $\chi$.
\end{theorem}
\begin{proof}
 Assume that $\chi(f)=\chi(g)$: then $f''/f'=g''/g'$,
 $g'(0)=f'(0)$ and $f(0)=g(0)=0$.  Thus the same reasoning as in
 the proof of Theorem \ref{th:hatbetaopen} implies that $f=g$.

 Since $f \in \mathcal{O}_{qc}$ implies $\alpha f \in \Oqc$
 for any $\alpha \in \mathbb{C}^*$
 and $\mathcal{A}(\alpha f)=\mathcal{A}(f)$,  if
 $(\psi,a) \in \chi(\Oqc)$ then so is $(\psi,\alpha a)$ for all
 $\alpha \in \mathbb{C}^*$.  Furthermore
 $\mathcal{A}(\Oqc)=\mathcal{A}(\widetilde{\mathcal{D}})=\hat{\beta}
 (\widetilde{\mathcal{D}})$.  It follows that
 $\chi(\Oqc)=\hat{\beta}(\widetilde{\mathcal{D}}) \oplus \mathbb{C}^*$.

 The last claim follows from the fact that
 $\hat{\beta}(\widetilde{\mathcal{D}})$ is open by Theorem
 \ref{th:hatbetaopen}.
\end{proof}

The metric on $A^1_\infty(\mathbb{D}) \oplus
\mathbb{C}$ induces the metric
\[ d_o(f,g)=||\chi(f)-\chi(g)||= || \mathcal{A}(f)- \mathcal{A}(g)
  ||_{1,\infty} + |f'(0)-g'(0)| \]
on $\Oqc$. With the topology induced by this metric we immediately
have
\begin{proposition} \label{pr:dividebyf'(0)} The map
\begin{eqnarray*}
 \Oqc & \rightarrow & \widetilde{\mathcal{D}} \\
 f & \mapsto & f/f'(0)
\end{eqnarray*}
is continuous.
\end{proposition}

 A few results about the topology of $\Oqc$
 will be necessary.  We will
 need the fact that point evaluation at $z$ is a continuous map on
 $\widetilde{\mathcal{D}}$ and $\Oqc$ for any $z \in \overline{\mathbb{D}}$.
 We will also need a theorem to the effect that if $f$ and $g$ are
 nearby in $\widetilde{\mathcal{D}}$ then $f(S^1)$ is
 uniformly close to $g(S^1)$, and similarly for $\Oqc$.

 In order to do this we will make use of the
 Teichm\"uller metric, which is defined as follows.  Let $f$ and
 $g$ be elements of $\mathcal{D}$, where $\mathcal{D}$
 is considered as a model for $T(\mathbb{D}^*)$.
 The Teichm\"uller distance
 between $f$ and $g$ is defined to be
 \[  \tau(f,g)= \frac{1}{2} \inf_{\mu,\nu} \left( \frac{1+||(\mu -
  \nu)/(1-\overline{\mu}\nu)||_\infty}{1 - ||(\mu -
  \nu)/(1-\overline{\mu}\nu)||_\infty} \right)   \]
 where $\mu$ and $\nu$ are the dilatations of quasiconformal
 extensions of $f$ and $g$ respectively, and the minimum is taken
 over all such choices of extension.  The Teichm\"uller metric is
 compatible with the Schwarzian metric $d_s$ introduced in Section
 \ref{se:BerspreBers} (see \cite[III.4.2, III.4.3]{Lehto}).
 Furthermore, if one fixes
 a quasiconformal extension $\tilde{f}$ of $f$ to
 $\overline{\mathbb{C}}$ with dilatation $\mu_0$, and takes the
 minimum over $\nu$ only, then the result is the same
 (see \cite[Lemma III.2.1]{Lehto}).
 That is,
 \begin{equation*}
  \tau(f,g)= \frac{1}{2} \inf_{\nu} \left( \frac{1+||(\mu_0 -
  \nu)/(1-\overline{\mu_0}\nu)||_\infty}{1 - ||(\mu_0 -
  \nu)/(1-\overline{\mu_0}\nu)||_\infty} \right).
 \end{equation*}
 Note that, if $K_{h}$ denotes the maximal dilatation of
 a quasiconformal map $h$, then
 \begin{equation}  \label{eq:Teichmullerdistanceonemu}
  \tau(f,g)= \frac{1}{2} \inf_{\tilde{g}} \log{K_{\tilde{g} \circ
  \tilde{f}^{-1}}}.
 \end{equation}
 where $\tilde{f}$ is a fixed quasiconformal extension of $f$ and
 $\tilde{g}$ are quasiconformal extensions of $g$ to
 $\overline{\mathbb{C}}$.
 \begin{theorem} \label{th:S1nearbyDtilde}
  The point evaluation map $f \mapsto f(z)$ is continuous on
  $\widetilde{\mathcal{D}}$ for any $z \in \overline{\mathbb{D}}$.
  Furthermore, for any $f \in \widetilde{\mathcal{D}}$ and $\epsilon>0$,
  there is an open neighbourhood
  $U \subset \widetilde{\mathcal{D}}$ containing $f$ such that for any
  $g \in U$ and $y \in \partial \mathbb{D}$, $|g(y)-f(y)| < \epsilon$.
 \end{theorem}
 \begin{proof}
  First, fix any $K>1$ and $R>1$.  Consider the family of K-quasiconformal
  mappings of the disc $|z|< R$ whose restriction to $\mathbb{D}$
  is in $\widetilde{\mathcal{D}}$.
  This family is equicontinuous in $|z|<R$ by \cite[Theorem 4.2]{Lehto-Virtanen},
  since the family omits $\infty$, and for any
  point $z$ such that $|z|<1$, $f(z)$ is uniformly bounded away from
  $\infty$ by the growth theorem for univalent mappings of the disc.

  Now fix $f \in \widetilde{\mathcal{D}}$.
  Fix a quasiconformal extension $\tilde{f}$ of $f$, and let
  $K>1$ be the maximal dilatation of $f$.  Finally
  fix a $K_1>K$.  We claim that there is a $\delta_1>0$ such that if
  $d(f,g)_{ps} < \delta_1$ then $g$ has an extension with
  maximal dilatation less than
  $K_1$.  To see this, choose M\"obius transformations $T_1$ and
  $T_2$ such that $T_1 \circ f,\  T_2 \circ g \in \mathcal{D}$.  By
  the compatibility of the Schwarzian metric $d_s$ on
  $\mathcal{D}$ and the Teichm\"uller metric $\tau$ in the form
  (\ref{eq:Teichmullerdistanceonemu}), there is a $\delta_2$ such
  that $d_s(T_1 \circ f,T_2 \circ g) <\delta_2$ implies that $T_2
  \circ g$ has a quasiconformal extension to
  $\overline{\mathbb{C}}$ with maximal dilatation less than or equal to
  $K_1$.  Since for any M\"obius transformation $T$,
  post-composition by $T$ does not affect the maximal dilatation of a
  quasiconformal map, and since $\mathcal{S}(T \circ f)=\mathcal{S}(f)$,
  it follows that if $||\mathcal{S}(f)-\mathcal{S}(g)||_{1,\infty} <
  \delta_2$ then $g$ has a quasiconformal extension of
  dilatation less than $K_1$.  Since $\Psi$ is continuous by
  Theorem \ref{th:Psicontinuous}, and $\Psi \circ
  \mathcal{A}=\mathcal{S}$, the claim follows.

  We now show that point evaluation is continuous in $\widetilde{\mathcal{D}}$
  for $z \in \overline{\mathbb{D}}$.  First observe that the claim is
  immediate for $z \in \mathbb{D}$ since convergence in $d_{ps}$ implies
  uniform convergence in compact subsets of $\mathbb{D}$.  To prove the
  claim for $x \in S^1$, fix $f \in \widetilde{\mathcal{D}}$
  and a quasiconformal extension $\tilde{f}$; this
  $\tilde{f}$ has maximal dilatation $K$ for some $K>1$.
  For $K_1>K$, recall that the family $\mathcal{F}$
  of mappings in $\widetilde{\mathcal{D}}$ with quasiconformal extensions
  to $|z|<R$ of maximal dilatation less
  than $K_1$ is equicontinuous in  $\overline{\mathbb{D}}$
  by the first paragraph of the proof. This family
  includes $f$ itself.  Furthermore, the boundary values of $g \in
  \mathcal{F}$ do not depend on the choice of extension.
  We show that point evaluation $g \mapsto g(x)$ is continuous at $f$.
  Choose a $y \in \mathbb{D}$ such that
  $|g(x)-g(y)|<\epsilon/3$ for any $g \in \mathcal{F}$; this can be
  done by the equicontinuity of $\mathcal{F}$.  There is a
  $\delta$ such that $d_{ps}(f,g) < \delta \Rightarrow |f(y)-g(y)|
  < \epsilon/3$.  Shrinking $\delta$ if necessary, by the previous
  paragraph we may take $\delta$ so that $d_{ps}(f,g) <\delta$
  implies that $g \in \mathcal{F}$.  Thus $|g(x)-f(x)|\leq |g(x)-g(y)|+
  |g(y)-f(y)|+|f(y)-f(x)| <\epsilon$.

  Since $\partial\mathbb{D}$ is compact, to
  prove the last claim of the Theorem
  it suffices to show that for any $x \in \partial\mathbb{D}$ and $\epsilon>0$, there is an open
  subinterval of $\partial\mathbb{D}$ containing $x$ and a $\delta$ such that
  $d_{ps}(f,g)<\delta \Rightarrow |g(y)-f(y)|<\epsilon$
  for all $y$ in the interval.
  By the equicontinuity of $\mathcal{F}$ there is an open subinterval $I$
  of $\partial\mathbb{D}$ containing $x$ and a $\delta$ such that for all $y \in I$, if $g$
  satisfies $d_{ps}(f,g)< \delta$ then $|g(y)-g(x)|<\epsilon/3$.
  Furthermore by the continuity of point evaluations, by shrinking
  $\delta$ if necessary we can ensure that
  $|g(x)-f(x)|<\epsilon/3$.  Thus another three-$\epsilon$
  argument
  \begin{equation*}
   |g(y)-f(y)| \leq |g(y)-g(x)| + |g(x)-f(x)|+|f(x)-f(y)|
   < \epsilon
  \end{equation*}
  completes the proof.
 \end{proof}

 \begin{corollary} \label{co:Oqcpointevalcontin}
 The point evaluation map $f \mapsto f(z)$ is continuous
 on $\Oqc$ for any $z \in \overline{\mathbb{D}}$.  Furthermore,
 for any $f \in \Oqc$ and $\epsilon >0$, there is
 an open neighbourhood $U \subset \Oqc$ of $f$ such that for any $g
 \in U$ and $y \in \partial\mathbb{D}$, $|g(y)-f(y)|< \epsilon$.
 \end{corollary}
\begin{proof}
 Point evaluation is continuous for $z \in \mathbb{D}$ since
 convergence in $\Oqc$ implies uniform convergence on compact
 subsets of $\mathbb{D}$.  The continuity of point evaluation on
 $\partial\mathbb{D}$ follows from the second claim of the Corollary, which we
 now prove.

 Fix $f \in \Oqc$.  By Proposition \ref{pr:dividebyf'(0)} and
 the second claim of Theorem \ref{th:S1nearbyDtilde}, there is a
 $\delta_1<0$ and $M>0$ such that if $d_o(f,g) < \delta_1$
 then $|g(y)|/|g'(0)|\leq M$ for all $y \in \partial\mathbb{D}$.
Applying Proposition \ref{pr:dividebyf'(0)} and
 Theorem \ref{th:S1nearbyDtilde},
 there is a $\delta_2$ such that $d_o(f,g)<\delta_2$ implies that
 \[  \left| \frac{g(y)}{g'(0)} - \frac{f(y)}{f'(0)}
   \right|<\frac{\epsilon}{2|f'(0)|} \]
 for all $y \in \partial\mathbb{D}$.  Since $\chi$ is continuous in
 its second component there is a $\delta_3$ such that
 $d_o(f,g)<\delta_3$ implies that
 \[  \left|\frac{g'(0)}{f'(0)} -1 \right|< \frac{\epsilon}{2M
    |f'(0)|}.  \]
 Choosing $\delta$ to be the minimum of $\delta_1$, $\delta_2$ and
 $\delta_3$, the claim follows from the estimate
 \begin{eqnarray*}
  |g(y)-f(y)| & = & |f'(0)| \left| \frac{g'(0)}{f'(0)}
  \frac{g(y)}{g'(0)} - \frac{g(y)}{g'(0)} + \frac{g(y)}{g'(0)} -
  \frac{f(y)}{f'(0)} \right| \\
  & \leq  & |f'(0)| \left( \left| \frac{g'(0)}{f'(0)} -1 \right|
  \left| \frac{g(y)}{g'(0)} \right| + \left| \frac{g(y)}{g'(0)} -
  \frac{f(y)}{f'(0)} \right| \right).
 \end{eqnarray*}
\end{proof}
Corollary \ref{co:Oqcpointevalcontin} immediately implies:
\begin{corollary} \label{co:nearnesslemma}
 Let $f \in \Oqc$ and $B$ be an open subset of
 $\overline{\mathbb{C}}$ containing $f(\overline{\mathbb{D}})$.
 There exists an open neighbourhood $U$ of $f$ in $\Oqc$ such that
 $g(\overline{\mathbb{D}}) \subset B$ for all $g \in U$.
\end{corollary}

\end{subsection}

\begin{subsection}{Complex structure on $\Oqc(\riem)$}
In this section we construct a complex Banach manifold structure on
$\Oqc(\riem)$.

First we outline some facts regarding infinite-dimensional
holomorphy (see for example  \cite{Chae}, \cite[V.5.1]{Lehto} or
\cite[Section 1.6]{Nag}). Let $E$ and $F$ be Banach spaces. Let $U$
be an open subset of $E$.
\begin{definition}  A map $f \colon U \rightarrow F$ is holomorphic if
 for each $x_0 \in U$ there is a continuous complex linear map
 $Df(x_0) \colon E \rightarrow F$ such that
 \[  \lim_{h \rightarrow 0} \frac{|| f(x_0+h)-f(x_0) -
 Df(x_0)(h)||_F}{||h||_E} =0.  \]
\end{definition}
\begin{definition}
A map $f  \colon U \to F$ is called G\^ateaux holomorphic if $f$ is
holomorphic on affine lines. That is, if for all $a \in U$ and all
$x \in E$, the map $z \mapsto  f(a+zx)$ is holomorphic on $\{z \in
\mathbb{C} \st a + zx \in U \}$.
\end{definition}

\begin{theorem}[{\cite[p 198]{Chae}}]
\label{th:holo_equivalents} Let $f  \colon U \to F$. The following are
equivalent.
\begin{enumerate}
\item $f$ is holomorphic.
\item $f$ is G\^{a}teaux-holomorphic and continuous.
\item $f$ is G\^{a}teaux-holomorphic and locally bounded on $U$.
\end{enumerate}
\end{theorem}
\begin{remark}  In \cite[Lemma 5.1 (i) p 206 ]{Lehto} it is stated that
 G\^ateaux holomorphy is sufficient for holomorphy.  In fact
 some further condition such as continuity or local boundedness
 is necessary.
\end{remark}

 In order to construct the complex structure on $\Oqc(\riem)$,
 we will need the following Lemma, which says that composition on
 the left by a holomorphic map is a holomorphic operation on
 small neighbourhoods in $\mathcal{O}_{qc}$.
 \begin{lemma} \label{le:chartlemma}
  Let $K$ be a compact set which is the closure of an open
  neighbourhood of $0$ in $\mathbb{C}$, and
  let $A$ be an open set containing $K$. Let $U$ be an open set in $\Oqc$ which satisfies
  $\psi(\overline{\mathbb{D}}) \subset K$ for all $\psi \in U$.  Let
  $h \colon A \rightarrow \mathbb{C}$ be a one-to-one holomorphic map
  satisfying $h(0)=0$.  The map
  \begin{align*}
   L_h \colon U & \rightarrow \Oqc \\
   f & \mapsto h \circ f
  \end{align*}
  is a holomorphic map from $U$ into $\Oqc$.
 \end{lemma}

The existence of a set $U$ satisfying the hypotheses of the Lemma is
established by Corollary \ref{co:nearnesslemma}.  We now proceed
with the proof.

 \begin{proof}
  Given the complex structure on $\Oqc$ provided by Theorem \ref{th:chi_induces_complex_structure},
  to show $L_h$ is holomorphic
  we must show that $\chi \circ L_h \circ \chi^{-1} \colon \chi(U) \to A^1_{\infty}
  \oplus \mathbb{C}$ is holomorphic. That is,
  $\left(\mathcal{A}(f), f'(0) \right) \mapsto \left( \mathcal{A}(h \circ f),
  (h \circ f)'(0) \right)$ is holomorphic. We prove separate
  holomorphicity and apply Hartogs' theorem (a suitable statement in infinite
  dimensions can be found in \cite{Mujica}).

  It is clear that the map $f'(0) \mapsto (h \circ f)'(0)=h'(0)
  f'(0)$ is holomorphic and continuous, so it remains to be proven that
  \[  \mathcal{A}(f) \mapsto \mathcal{A}(h \circ f) \]
  is holomorphic in
  $A^1_\infty(\mathbb{D})$.  The map is locally (in fact globally) bounded, since
  $\mathcal{A}(\Oqc)$ is a bounded subset of
  $A^1_\infty(\mathbb{D})$ by the elementary inequality
  \[  (1-|z|^2)\left|\frac{f''(z)}{f'(z)}\right| \leq 6  \]
  for univalent functions $f$ satisfying $f(0)=0$.
  Thus by Theorem \ref{th:holo_equivalents} we need
  only show that this map is G\^ateaux holomorphic.

  Choose $f_0 \in U$.  We need to show that
  $t \mapsto \mathcal{A}(h \circ f_t)$ is a holomorphic curve in
  $A^1_\infty(\mathbb{D})$ for complex curves $f_t \in
  \Oqc$ of the following form: $f_t$ is the
  inverse image under $\chi$ of $(\mathcal{A}(f_0) + t \phi,q(t))$ where $\phi
  \in A^1_\infty(\mathbb{D})$ and $q$ is holomorphic in $t$ with $q(0)=f_0'(0)$.
  Since $\chi$ is continuous and $\chi(\Oqc)$ is open by Theorem
  \ref{th:chi_induces_complex_structure}, it follows that $f_t \in
  U$ for $t$ in a small enough neighbourhood of $0$.
  The explicit solution to this
  differential equation can be
  obtained in terms of $\phi$, $f_0$ and $q(t)$:
  \[  f_t(z)=\frac{q(t)}{f_0'(0)} \int_0^z f_0'(u) \exp{\left( t\int
       _0^u \phi(w)dw \right)} du.  \]
  Note that for fixed $z$, $f_t(z)$ is holomorphic in $t$.
  Denoting differentiation with respect to $t$ by $\dot{f_t}$, this
  expression also implies that $\dot{f}_t$, $\dot{f}'_t$,
  $\ddot{f}_t$ and $\ddot{f}'_t$ exist and are holomorphic in $z$.

  Let $f_0 \in U$.  Assume that we have a curve $f_t$ of the above form
  for $t \in B$, where $B\subset \mathbb{C}$ is some open neighbourhood
  of $0$.  Shrinking $B$
  if necessary,  Theorem \ref{th:S1nearbyDtilde} guarantees that
  $f_t(\overline{\mathbb{D}}) \subset K$ for all $t \in B$; that is,
  $f_t \in U$.  Define
  \[  \alpha(t) = \frac{h''}{h'} \circ f_t \cdot f_t'.  \]
  Denoting $t$-differentiation with a dot we then have that
  \[  \lim_{t \rightarrow 0} \frac{1}{t} \left( \mathcal{A}(h \circ
  f_t) - \mathcal{A}(h \circ f_0) \right) = \dot{\alpha}(t) + \phi
  \]
  where
  \[  \dot{\alpha}(t) = \left( \frac{h'''}{h'}-
  \left(\frac{{h''}}{h'}\right)^2
  \right) \circ f_t \cdot \dot{f}_t\, f_t' + \frac{h''}{h'} \circ f_t
  \cdot \dot{f}_t'.  \]

  It is easy to compute that
  \[  \left\|\frac{1}{t} \left( \mathcal{A}(h \circ f_t) - \mathcal{A}(h
    \circ f_0)
   \right) -  (\dot{\alpha}(t) + \phi)  \right\|_{1,\infty} = \left\|
   \frac{1}{t}\left( \alpha(t) - \alpha(0) - t \dot{\alpha}(t) \right)
   \right\|_{1,\infty},  \]
  since $\mathcal{A}(f_t) - \mathcal{A}(f_0)=t \phi$.
  We will show that the right
  hand side goes to $0$ as $t \rightarrow 0$.

  By an elementary estimate
  \[  |\alpha(t) - \alpha(0) - t \dot{\alpha}(t)| \leq \sup_{|s| \leq
      |t|} | \ddot{\alpha}(s) | t^2,  \]
  so it only need be shown that $\ddot{\alpha}(t) \in
  A^1_\infty(\mathbb{D})$ for $t$ in some neighbourhood of $0 \in
  \mathbb{C}$.  It is easily computed that
  \begin{eqnarray} \label{eq:alphadoubledot}
   \ddot{\alpha}(t) & = & \left( \frac{h''''}{h'} - 3 \frac{h'''
   h''}{{h'}^2} - \frac{{h''}^3}{{h'}^3} \right) \circ f_t \cdot
   {\dot{f}_t}^2 \, f_t' + \left( \frac{h'''}{h'} -
   \frac{{h''}^2}{{h'}^2} \right)
   \circ f_t \cdot  \left( \ddot{f}_t \, f_t' + 2 \dot{f}_t \,
    {\dot{f}_t}' \right) \\
    &  & + \frac{h''}{h'} \circ f_t \cdot \ddot{f}_t' \nonumber
  \end{eqnarray}

  We claim that there exists an $r>0$ and positive constants
  $M_1$, $M_2$, $M_3$ which bound
  $f_s(z)$, $\dot{f}_s(z)$ and $\ddot{f}_s(z)$ respectively
  for all $z \in \mathbb{D}$ and
  $|s| < r$.  We know by Theorem \ref{th:S1nearbyDtilde}
  that there is a constant $M_1$ and $r_1>0$ such that $|f_s(z)|
  \leq M_1$ for $|s|<r_1$ and $z \in \mathbb{D}$.
  Let $r_2$ and $r'$ satisfy $0< r_2 < r' < r_1$. Applying Cauchy
  estimates in the variable $s$ using the curve $|s|=r'$ we have
  that $|\dot{f}_s(z) | \leq r' |f_s(z)|/|r_1-r'|^2 \leq M_2$ for
  some $M_2$ which is independent of $z \in \mathbb{D}$ and $|s|<r_2$.  The same
  argument again shows that there is an $r>0$ and an $M_3$ such that
  $|\ddot{f}_s(z)| \leq M_3$ for $|s|<r$ uniformly in $z$,
  so the claim follows.

  Applying the Schwarz lemma to $f_t$, we have that $(1-|z|^2)
  |f_t'(z)| \leq (M_1^2-|f(z)|^2)/M_1 \leq M_1$ for $|t|<r$, and similarly
  $(1-|z|^2)
  |\dot{f}_t'(z)| \leq M_2$ and $(1-|z|^2) |\ddot{f}_t'(z)| \leq M_3$ for
  all $|t|<r$ and $z \in \mathbb{D}$.  Since $h$ is holomorphic on
  $K$  and $h' \neq 0$, the
  factors $(h''/h') \circ f_t$ etc. of each term in equation (\ref{eq:alphadoubledot})
  are bounded by constants
  $C_1$, $C_2$ and $C_3$ respectively, uniformly for $z \in \mathbb{D}$ and $|t|<r$.
  Applying these estimates to equation (\ref{eq:alphadoubledot})
  we have
  \[  (1-|z|^2)|\ddot{\alpha}(s)(z)| \leq C_1 M_1^2M_2 + C_2
  (M_3M_1+2M_2^2)  + C_3 M_3  \]
  which completes the proof.
 \end{proof}

 We may now state and prove the main theorem.
 \begin{theorem} \label{th:Oqcriem_complexstructure} Let
 $\riem$ be a Riemann surface with $n$ distinguished points $p_1, \ldots, p_n$.
 The space $\Oqc(\riem)$ has a complex Banach manifold
 structure modelled on $\Oqc^n=\Oqc \times \cdots \times \Oqc$.
 \end{theorem}
 \begin{proof}
 Fix a point $(\phi_1,\ldots,\phi_n) \in \Oqc(\riem)$.  We will
 construct an open neighbourhood and a chart in $\mathcal{O}_{qc}(\riem)$.

 For $i = 1,\dots,n$, let $D_i=\phi_i(\mathbb{D})$ be open sets in $\riem$.  Choose
 domains $B_i \subset \riem$ with the following properties: 1)
 $\phi_i(\overline{\mathbb{D}}) \subset B_i$, 2) $B_i \cap B_j =
 \emptyset$ for $i \neq j$ and 3) $B_i$ are open and simply connected.
 Let $\zeta_i \colon B_i \rightarrow \mathbb{C}$ be a local biholomorphic
 parameter such that $\zeta_i(p_i)=0$ and $\zeta_i(B_i)$ is bounded.  We then have that $\zeta_i
 \circ \phi_i \in \Oqc$.

 Let $K_i$ be a compact set containing $\zeta_i \circ
 \phi_i(\overline{\mathbb{D}})$ in its interior,
 such that $\zeta_i^{-1}(K_i) \subset B_i$.
 By Corollary \ref{co:nearnesslemma}, for each $i$
 there is an open neighbourhood $U_i$ of $\zeta_i \circ \phi_i$
 in $\Oqc$ such that $\psi_i(\overline{\mathbb{D}}) \subset K_i$ for
 all $\psi_i \in U_i$, and so $(\zeta_1^{-1} \circ \psi_1, \ldots,
 \zeta_n^{-1} \circ \psi_n)$ is an element of $\Oqc(\riem)$.

Thus, $U_1 \times \cdots \times U_n$ is an open subset of $\Oqc
\times \cdots \times \Oqc$ with the product topology.  Let
\[ V_i=\{ \,\zeta_i^{-1} \circ \psi_i \,|\, \psi_i \in U_i \,\},  \]
and
\begin{align*}
 T_i \colon V_i  & \rightarrow \Oqc \\
 g & \mapsto \zeta_i \circ g.
\end{align*}
Of course we have that $T_i(V_i)=U_i$.

Let $\Oqc(\riem)$ be endowed with the topology generated by the
open sets of the above form. We define the complex structure on
$\Oqc(\riem)$ as follows. Let $V=V_1 \times \cdots \times V_n
\subset \Oqc(\riem)$.  The charts are ordered pairs $(T,V)$ where
$T=T_1 \times \cdots \times T_n \colon V \rightarrow \Oqc^n$.
Consider two charts $(T,V)$ and $(T',V')$ such that $V \cap V'
\neq \emptyset$. Let $\zeta_i$ and $\zeta_i'$, $K_i$ and $K_i'$,
etc. correspond to the charts $(T,V)$ and $(T',V')$ respectively.
We have that $\zeta_i' \circ \zeta_i^{-1}$ is holomorphic on
$K_i$, and $T_i' \circ T_i^{-1}(\psi_i) = \zeta_i' \circ
\zeta_i^{-1} \circ \psi_i$ for any $\psi_i \in T_i(V_i \cap V_i')
\subset \Oqc$. By Lemma \ref{le:chartlemma} the map $\psi_i
\mapsto \zeta' \circ \zeta^{-1} \circ \psi_i$ is holomorphic on
$T_i(V_i \cap V_i')$. Applying the argument to $T_i \circ
T_i'^{-1}$ shows that it is a biholomorphism. Thus $T' \circ
T^{-1}$ is a biholomorphism.
\end{proof}

\begin{remark} The topology on $\Oqc(\riem)$ defined in the proof of
Theorem \ref{th:Oqcriem_complexstructure} is Hausdorff. This can be
proved with a short direct argument that relies on the following two
facts. (1) $\Oqc$ is Hausdorff. (2) If $(\phi_1,\ldots,\phi_n)$ and $(\phi'_1,\ldots,\phi'_n)$ are in
$\Oqc(\riem)$ then biholomorphic parameters $\zeta_i = \zeta'_i$ can be chosen
such that the neighbourhoods $B_i$ and $B'_i$ of the images of $\phi_i$
and $\phi'_i$ are contained in the domain of $\zeta_i$.
\end{remark}

\end{subsection}
\end{section}


\begin{section}{Appendix: Identification of $\Oqc$ with $QS(S^1)
\times \mathbb{R}$} \label{se:Oqcappendix}
 We describe the relation of $\Oqc$ to $\mbox{QS}(S^1)$.

In the following, we say that two domains $E$, $F$ in
$\overline{\mathbb{C}}$ are complementary if $F$ is the complement
of the closure of $E$.  In particular $\mathbb{D}^*$ and
$\mathbb{D}$ are complementary.
\begin{theorem}[conformal welding] \label{conformalweldingone}
 Let $\gamma \in \mbox{QS}(S^1)$.  There is a unique $f \colon\mathbb{D}
 \rightarrow \mathbb{C}$ and $g \colon\mathbb{D}^* \rightarrow
 \mathbb{\overline{C}}$ which map onto complementary domains
 satisfying
 \begin{enumerate}
 \item $f$ and $g$ are one-to-one, holomorphic, and have quasiconformal extensions to
   $\overline{\mathbb{C}}$ 
\item $\gamma=\left. g^{-1} \circ f\right|_{S^1}$
 \item $f(0)=0$, $|f'(0)|=1$ 
 \item $g(\infty)=\infty$, $g'(\infty) >0$.
 \end{enumerate}
\end{theorem}
\begin{proof}
 This is a standard result with a different choice of
 normalization.  By \cite[Theorem 2.3]{Teo}
 there are maps $f$ and $g$ such that $(1)$ and $(4)$ hold but
 with $(2)$ and $(3)$ replaced by $(2')\  \gamma = \left.g^{-1} \circ f\right|_{S^1}
 \mbox{mod}(S^1)$ and $(3')\  f(0)=0, \ f'(0)=1$.  Property $(2')$ means
 that there
 exists some $\alpha \in \mbox{M\"ob}(S^1)$ corresponding to a rotation of
 $S^1$ such that $\gamma= \left. \alpha \circ g^{-1} \circ
 f\right|_{S^1}$.  Since $f$ and $g$ are unique, $\alpha$ is also
 uniquely determined.  Setting $\tilde{f}=\alpha^{-1} \circ f$ and
 $\tilde{g} = \alpha^{-1} \circ g \circ \alpha$ we have that
 $\tilde{f}$ and $\tilde{g}$ are the unique maps satisfying
 properties $(1)$ to $(4)$.
\end{proof}

\begin{corollary} \label{conformalweldingm}
 Let $m \in \mathbb{R}$ and $\gamma \in \mbox{QS}(S^1)$.  There is a
 unique $f \colon\mathbb{D}\rightarrow \mathbb{C}$ and $g \colon\mathbb{D}^*
 \rightarrow \overline{\mathbb{C}}$ mapping onto complementary
 domains and satisfying
 \begin{enumerate}
 \item $f$ and $g$ are one-to-one, holomorphic, and have
 quasiconformal extensions to
 $\overline{\mathbb{C}}$

  \item $\gamma = \left. g^{-1} \circ f\right|_{S^1}$

  \item  $f(0)=0$, $|f'(0)|=e^m$

 \item  $g(\infty)=\infty$, $g'(\infty) >0$.
 \end{enumerate}
\end{corollary}
\begin{proof}
 Apply Theorem \ref{conformalweldingone} to get $\tilde{f}$ and
 $\tilde{g}$ satisfying its conditions.  Then $f=e^m\tilde{f}$ and
 $g=e^m \tilde{g}$ satisfies the desired conditions.
\end{proof}
Corollary \ref{conformalweldingm} establishes the existence of the
bijection between $QS(S^1) \times \mathbb{R}$ and $\Oqc$.
\begin{align*}
 \mathcal{V} \colon \mbox{QS}(S^1) \times \mathbb{R} & \rightarrow \Oqc \\
 (\gamma,m) & \mapsto f \\
\end{align*}
The inverse of this mapping is
\[   \mathcal{V}^{-1}(f) = (\left. g^{-1} \circ f
\right|_{S^1},\log{|f'(0)|}). \]

\end{section}

\end{document}